# PARAMETER ESTIMATES FOR FRACTIONAL AUTOREGRESSIVE SPATIAL PROCESSES

By Y. Boissy, B. B. Bhattacharyya, X. Li and G. D. Richardson

*University of Central Florida, North Carolina State University, University of Central Florida and University of Central Florida*

A binomial-type operator on a stationary Gaussian process is introduced in order to model long memory in the spatial context. Consistent estimators of model parameters are demonstrated. In particular, it is shown that $\hat{d}_N - d = O_P(\frac{(\text{Log } N)^3}{N})$, where $d = (d_1, d_2)$ denotes the long memory parameter.

**1. Introduction and main results.** A process obeying the spatial autoregressive model

$$(1.1) \qquad X_{st} = \alpha X_{s-1,t} + \beta X_{s,t-1} - \alpha\beta X_{s-1,t-1} + \varepsilon_{st}$$

was first studied by Martin [18], where $-1 < \alpha, \beta < 1$. Martin indicated that it is often desirable in practice that a process be reflection symmetric, that is, $\rho_{k\ell} = \rho_{-k,-\ell} = \rho_{k,-\ell} = \rho_{-k,\ell}$ for lags $k$ and $\ell$, and that the autocorrelations have a simple form. These requirements led to the definition of model (1.1), which has autocorrelation $\rho_{k,\ell} = \alpha^{|k|}\beta^{|\ell|}$ for lags $k$ and $\ell$. Spatial autoregressive models are shown by Tjøstheim [22] to be useful in studying geophysical quantities such as seismological data. Jain [14] indicates that these models can be applied to develop useful algorithms for image processing. Culles and Gleeson [8] and Basu and Reinsel [2] illustrate the suitability of model (1.1) as an error term in a regression model used to analyze data collected in agricultural field trials. Moreover, empirical evidence of slow decay of correlations between yield in two-dimensional agricultural field trials has received considerable attention (e.g., [10, 24, 25] and [20]). This led to the study of power law correlation functions by Whittle [25] and Besag [3] as an alternative to exponential decay. Martin [19] also indicates the importance of long-range correlation structures in agricultural field experiments.

In addition, Professor Kwang-Yul Kim of the Department of Meteorology at Florida State University has communicated to us that many geophysical variables, such as ocean temperature, exhibit a well-extended spatial correlation structure (e.g., [16, 17]). The above references motivate the possible need for inclusion of a long memory component when modeling certain spatial processes.

It is assumed that the data-sites can be arranged in a square lattice. This commonly occurs in agricultural experiments where exactly one plant is located at site $(i, j)$. Tjøstheim [21] mentions that irregularly spaced data-sites can sometimes be replaced by a regular grid using the interpolation techniques of Delfiner and Delhomme [9].

Let $\mathbb{C}$ denote the field of complex numbers; define

$$
\begin{aligned}
\phi(z_1, z_2, \alpha, \beta) &= (1 - \alpha z_1)(1 - \beta z_2), \\
\psi(z_1, z_2, \alpha, \beta) &= [\phi(z_1, z_2, \alpha, \beta)]^{-1}.
\end{aligned}
\tag{1.2}
$$

Suppose that $B_1 X_{st} = X_{s-1,t} (B_2 X_{st} = X_{s,t-1})$ is the backward shift operator on the first (second) index of $X_{st}$; then (1.1) can be written compactly as

$$
\phi(B_1, B_2, \alpha, \beta) X_{st} = \varepsilon_{st}. \tag{1.3}
$$

It is assumed throughout this work that all processes considered are Gaussian. Asymptotic results on parameter estimators in the model (1.3) for both stationary and nonstationary cases can be found in [1, 4, 5, 15].

Our purpose here is to extend the work of Fox and Taqqu [11, 12] from time series to the spatial context by including a long memory component in the model (1.3). Memory in time series is modeled by use of the binomial operator $(1 - B)^d$, where $d \in (-\frac{1}{2}, \frac{1}{2})$ denotes the memory parameter. The following operator is used in the spatial setting with two indices $d = (d_1, d_2)$:

$$
\nabla^d = (1 - B_1)^{d_1} \circ (1 - B_2)^{d_2}. \tag{1.4}
$$

The operator $\nabla^d$ is defined by its corresponding power series representation in $(z_1, z_2)$, that is, $\nabla^d X_{st} = \sum_{k,\ell=0}^{\infty} a_{k\ell} X_{s-k, t-\ell}$, where the coefficients are found from the power series expansion of $(1 - z_1)^{d_1}(1 - z_2)^{d_2}$ in the unit polydisc $\Delta_1(0) \times \Delta_1(0)$ with $\Delta_1(0) = \{z \in \mathbb{C} : |z| < 1\}$. Let $\mathbb{Z}$ denote the set of all integers. Given $\phi$ in (1.2), $\nabla^d$ in (1.4) and the white noise process $\varepsilon_{st} \sim \mathcal{WN}(0, \sigma^2)$, the fractional autoregressive model of the form

$$
(1.5) \ \phi(B_1, B_2, \alpha, \beta) \nabla^d X_{st} = \varepsilon_{st}, \qquad \text{where } s, t \in \mathbb{Z} \text{ and } d_1, d_2 \in (-\tfrac{1}{2}, \tfrac{1}{2}),
$$

is considered. It can be shown (see [6]) that the spectral density function of the stationary solution of (1.5) is

$$
f(x, y, \theta) = \frac{\sigma^2}{4\pi^2} \frac{|1 - e^{-ix}|^{-2d_1} \cdot |1 - e^{-iy}|^{-2d_2}}{|\phi(e^{-ix}, e^{-iy}, \alpha, \beta)|^2}, \tag{1.6}
$$

where $\theta = (\alpha, \beta, d)$, $|\alpha| < 1, |\beta| < 1$ and $|d_i| < \frac{1}{2}$, $i = 1, 2$. Moreover, the corresponding autocovariance function when $\alpha = \beta = 0$ is

$$\gamma(k, \ell) = \frac{(-1)^{k+\ell}\Gamma(1 - 2d_1)\Gamma(1 - 2d_2)\sigma^2}{\Gamma(k - d_1 + 1)\Gamma(1 - k - d_1)\Gamma(\ell - d_2 + 1)\Gamma(1 - \ell - d_2)}.$$

Let $\theta_0 = (\alpha_0, \beta_0, d_0)$ denote the true parameter value, where $d_0 = (d_{10}, d_{20})$.

Recall that a second-order process $\{X_{st} : s, t \in \mathbb{Z}\}$ is said to be *stationary* when $E(X_{st}) = \mu$ and, for each $k, \ell \in \mathbb{Z}$, $\operatorname{cov}(X_{s+k, t+\ell}, X_{st})$ is independent of $s, t \in \mathbb{Z}$. Moreover, $\{X_{st} : s, t \in \mathbb{Z}\}$ is called *strictly stationary* provided all its finite-dimensional distributions remain invariant under translations. Denote $\overline{X}_N = \frac{1}{N^2}\sum_{k,\ell=1}^{N} X_{k\ell}$ and define

$$I_N(x, y) = \frac{1}{N^2} \left| \sum_{k,\ell=1}^{N} e^{i(kx+\ell y)} \cdot (X_{k\ell} - \overline{X}_N) \right|^2$$

to be the *periodogram* of the process. As suggested by Whittle [23] in the time series case, define $\hat{\theta}_N$ to be an argument $\theta = (\alpha, \beta, d)$ making

$$(1.7) \qquad \sigma_N^2(\theta) = \int_{I^2} \frac{I_N(x, y)}{f(x, y, \theta)} \, dx \, dy$$

a minimum, where $I = [-\pi, \pi]$, $I^2 = [-\pi, \pi] \times [-\pi, \pi]$ and $f(x, y, \theta)$ is the spectral density function of the process. In our setting, it can be shown that a minimum exists when $|\alpha| \leq r$, $|\beta| \leq r$ and $|d_i| \leq s$, where $0 < r < 1$ and $0 < s < \frac{1}{2}$, $i = 1, 2$.

Some motivational comments concerning the estimator defined in (1.7) are in order. Assume that $\{X_{st} : s, t \in \mathbb{Z}\}$ is a mean zero Gaussian process with spectral density function $f(x, y, \theta)$ given in (1.6) and denote $\mathbf{X} = (X_{11}, X_{12}, \ldots, X_{NN}) \sim \mathcal{N}(\mathbf{0}, \Gamma_N)$. Let $\operatorname{Log} x$ denote the natural logarithm of $x$. Then the log-likelihood function of $\mathbf{X}$ is

$$L_N(\theta, \sigma^2) = -\frac{N^2}{2}\operatorname{Log} 2\pi - \frac{1}{2}\operatorname{Log}|\Gamma_N| - \frac{1}{2}\mathbf{X}'\Gamma_N^{-1}\mathbf{X}.$$

Let $V_N = \sigma^{-2}\Gamma_N$. Whittle [23] proves that $|V_N(\theta)| \to 1$ as $N \to \infty$ and, thus, $L_N(\theta, \sigma^2)$ can be approximated by $-\frac{N^2}{2}\operatorname{Log} 2\pi\sigma^2 - \frac{1}{2\sigma^2}\mathbf{X}'V_N^{-1}(\theta)\mathbf{X}$ when $N$ is large. Fox and Taqqu ([11], page 518) use Parseval's identity to show that $\Gamma_N^{-1}$ can be approximated by the more tractable Toeplitz matrix $B_N = (b_{k\ell mn})$, where $b_{k\ell mn} = \int_{I^2} e^{i[(k-m)x+(\ell-n)y]} f^{-1}(x, y, \theta) \, dx \, dy$, and this leads to the Whittle estimator given in (1.7).

The primary results of this work are given below and proved in Sections 2–4. The many details and statements whose proofs are not given here are all available on request from the authors. Our first result establishes strong consistency of the estimator $\hat{\theta}_N$.

THEOREM 1.1. *Assume that $\{X_{st}: s, t \in \mathbb{Z}\}$ is a stationary Gaussian process having the spectral density function listed in* (1.6), *where* $|\alpha| < 1$, $|\beta| < 1$ *and* $0 < d_1, d_2 < \frac{1}{2}$. *Then* $\hat{\theta}_N \to \theta_0$ *almost surely as* $N \to \infty$.

The rate at which $\{\hat{\theta}_N\}$ converges in probability is given in the next theorem. The notation $\frac{\partial}{\partial \theta}\sigma_N^2(\theta_0)$ is defined to be the partial derivative of $\sigma_N^2$ evaluated at $\theta_0$.

THEOREM 1.2. *Let $\{X_{st}: s, t \in \mathbb{Z}\}$ be a stationary Gaussian process having the spectral density function listed in* (1.6), *where* $|\alpha| < 1$, $|\beta| < 1$ *and* $0 < d_1, d_2 < \frac{1}{2}$. *Then:*

(a) $N(\hat{\theta}_N - \theta_0) + NA_N^{-1}E(\frac{\partial}{\partial\theta}\sigma_N^2(\theta_0)) \xrightarrow{D} \mathcal{N}(\underline{0}, 128\pi^6\Sigma^{-1})$, *where $\Sigma$ and $A_N$ are defined in* (3.1), (3.2) *and* (3.5);

(b) $N(\operatorname{Log} N)^{-3}(\hat{\theta}_N - \theta_0) \xrightarrow{P} \frac{64\pi^2 K}{3}(\delta_{13}d_{10}, \delta_{24}d_{20}, \delta_{33}d_{10}, \delta_{44}d_{20})$, *where $K$ and $\delta_{ij}$ are defined in* (3.3).

Observe that $A_N$ in Theorem 1.2(a) depends on the abstract point $\overline{\theta}_N$ defined in (3.4). It follows that Theorem 1.2(a) cannot be used for construction of a confidence set about $\theta_0$. It only shows that when $\hat{\theta}_N - \theta_0$ is normalized by $N$, the second term in (a) provides just enough cancellation to produce a nondegenerate limit distribution.

REMARK 1.3. Many of the results given here are natural extensions of those given in the time series context by Fox and Taqqu [11]. However, a significant departure occurs with the bias term used in Theorem 1.2(a). According to (4.3), the normalized convergence of the sequence $\sum_{|k|<N}(1-|k|/N)a_{k1}c_{k1}$ is required. The proof given by Fox and Taqqu ([11], page 529) shows that (in our notation) $N^\rho \sum_{|k|<N}(1-|k|/N)a_{k1}c_{k1} \to 0$ and $N \to \infty$, for $\rho < 1$. (They only used the case $\rho = 1/2$.) In the spatial case one has to deal with the asymptotic behavior of $NE(\frac{\partial}{\partial\theta}\sigma_N^2(\theta_0))$, which amounts to considering what happens to $N\sum_{|k|<N}(1-|k|/N)a_{k1}c_{k1}$ as $N \to \infty$. Unlike in the time series setting [where the bias term $\sqrt{N}\sum_{|k|<N}(1-|k|/N)a_{k1}c_{k1} \to 0$], this sequence goes to infinity. It is shown in this paper that $(N/\operatorname{Log}^3 N) \times \sum_{|k|<N}(1-|k|/N)a_{k1}c_{k1}$ converges to a nonzero real number as $N \to \infty$. The point here is that the bias term, $NA_N^{-1}E[\frac{\partial}{\partial\theta}\sigma_N^2(\theta_0)]$, in the spatial context is of order $\operatorname{Log}^3 N$ as $N \to \infty$.

**2. Proof of Theorem 1.1: Outline.** A sequence of lemmas needed to establish the strong consistency of estimators of model parameters is stated and the reader is referred to Boissy [6] for detailed proofs. The first lemma gives a convenient method (in our setting) for finding the Fourier series representation of a spectral density function in terms of the orthogonal set $\{e^{i(kx+\ell y)} : k, \ell \in \mathbb{Z}\}$ with respect to the product Lebesgue measure $\lambda = \lambda_1 \times \lambda_2$ on $I^2$. Let $L^p(I^2)$ denote the set of all complex-valued functions for which $\int_{I^2} |f|^p d\lambda < \infty$, where $p > 0$.

LEMMA 2.1. *Let $g$ be a complex-valued function defined on the closed unit polydisc $\overline{\Delta_1(0)} \times \overline{\Delta_1(0)}$ in $\mathbb{C} \times \mathbb{C}$ for which $g(e^{-ix}, e^{-iy}) \in L^2(I^2)$. Suppose that $g$ is an analytic function defined on $\Delta_1(0) \times \Delta_1(0)$ with power series representation $g(z_1, z_2) = \sum_{k,\ell=0}^{\infty} C_{k\ell} z_1^k z_2^{\ell}$, and assume that $\lambda\{(x,y) \in I^2 : g$ is discontinuous at $(z_1, z_2) = (e^{-ix}, e^{-iy})\} = 0$. Then $g(e^{-ix}, e^{-iy})$ has the Fourier series expansion $\sum_{k,\ell=0}^{\infty} C_{k\ell} e^{-i(kx+\ell y)}$ if and only if $\sum_{k,\ell=0}^{\infty} |C_{k\ell}|^2 < \infty$.*

The representation described in Lemma 2.1 can be used to verify Lemma 2.2(a).

LEMMA 2.2. *Let $f(x, y, \theta)$ denote the spectral density function listed in (1.6). Then:*

(a) *$f \in L^1(I^2)$ and $\int_{I^2} \text{Log} \frac{4\pi^2}{\sigma^2} f(x, y, \theta) \, dx \, dy = 0$,*
(b) *$\int_{I^2} \frac{f(x,y,\theta_1)}{f(x,y,\theta_2)} \, dx \, dy > 4\pi^2$ when $\theta_1 \neq \theta_2$.*

LEMMA 2.3. *Assume that $\{X_{st} : s, t \in \mathbb{Z}\}$ is a stationary process defined on the underlying probability space $(\Omega, \mathfrak{F}, P)$ and having mean $\mu$, autocovariance $\gamma$ and spectral density function $h(x, y, \theta_0)$. Suppose that $\overline{X}_N \to \mu$ and $\frac{1}{N^2} \sum_{s,t=1}^{N} (X_{s+k,t+\ell} - \mu)(X_{st} - \mu) \to \gamma(k, \ell)$ almost surely, and that $g : (\mathbb{R}^{p+2}, \mathfrak{B}_{p+2}) \to (\mathbb{R}, \mathfrak{B})$ is bounded, Borel measurable and periodic with $g(\underline{x} + 2\pi \underline{k}) = g(\underline{x})$ for each $\underline{x} \in \mathbb{R}^{p+2}$ and $\underline{k} \in \mathbb{Z}^{p+2}$. If $I^2 \times K$ is a compact subset of $I^{p+2}$ which is contained in the set of all continuity points of $g$ and $I_N$ is the periodogram of the process, then $\int_{I^2} g(x, y, \theta) I_N(x, y) \, dx \, dy \to \int_{I^2} g(x, y, \theta) h(x, y, \theta_0) \, dx \, dy$ uniformly in $\theta \in K$, almost surely $[P]$.*

The omitted proof uses the technique of Hannan ([13], Lemma 1) in the time series case by uniformly approximating $g(x, y, \theta)$ over $I^2 \times K$ with an $N$th order Cesàro sum (see [26], page 304, Theorem 1.20). Since a stationary Gaussian process is strictly stationary, it can be shown that $\overline{X}_N \to \mu$ and $\frac{1}{N^2} \sum_{s,t=1}^{N} (X_{s+k,t+\ell} - \mu)(X_{st} - \mu) \to \gamma(k, \ell)$ almost surely and, thus, Lemma 2.3 is applicable in our context.

PROOF OF THEOREM 1.1. Combining Lemmas 2.2–2.3 above with the argument given by Hannan ([13], Theorem 1) for the time series case shows that $\hat{\theta}_N \to \theta_0$ almost surely (see [6]). □

Next, a sequence of estimators of $\theta_0$ that is more suited for computational purposes than (1.7) is given. Following Hannan ([13], page 133), define $\tilde{\theta}_N$ to be the argument making

$$(2.1) \qquad \tilde{\sigma}_N^2(\theta) = \frac{1}{N^2} \sum_{-N/2 < s,t \leq N/2} \frac{I_N(w_s, w_t)}{(4\pi^2/\sigma^2) f(w_s, w_t, \theta)}$$

a minimum, where $w_s = 2\pi s/N$. Under the hypothesis of Theorem 1.1 above, it can be shown that $\tilde{\theta}_N \to \theta_0$ almost surely. Once $\tilde{\theta}_N$ has been found, $\tilde{\sigma}_N^2(\tilde{\theta}_N)$ in (2.1) can be used as a strongly consistent estimator of $\sigma^2$.

**3. Proof of Theorem 1.2(a).** In the time series context, Fox and Taqqu [12] used a combinatorial method in order to establish convergence to normality of certain sequences of quadratic forms determined by Toeplitz matrices. Lemma 3.2 below allows one to extend the above result to the spatial context when the coefficient matrix can be expressed as a finite sum of Kronecker products. This lemma is used to prove Theorem 1.2(a).

Given $h \in L^1(I^2)$, $a_{k\ell mn} = \int_{I^2} e^{i[(k-m)x+(\ell-n)y]} h(x,y)\,dx\,dy$ is a Fourier coefficient of $h$ and $T_N(h) = (a_{k\ell mn})$ denotes the corresponding $N^2 \times N^2$ (block) Toeplitz matrix. Elements in $(a_{k\ell mn})$ are arranged lexicographically beginning with row 1, followed by row 2, and so on. In particular, element $a_{k\ell mn}$ appears in row $(k-1)N + \ell$ and column $(m-1)N + n$ of $T_N(h)$. Recall that the *Kronecker product* of matrices $A$ and $B$ is defined by $A \otimes B = (a_{k\ell}B)$, and the following properties are needed: $\prod_{j=1}^k (A_j \otimes B_j) = (\prod_{j=1}^k A_j) \otimes (\prod_{j=1}^k B_j)$ and $\mathrm{Tr}(A \otimes B) = \mathrm{Tr}(A) \cdot \mathrm{Tr}(B)$, when the matrices are compatible and Tr denotes the trace.

DEFINITION 3.1. A function $h: I \to \mathbb{R}$ is called *admissible* provided the following conditions are fulfilled:

(A.1) $h$ is symmetric and integrable;
(A.2) the set of discontinuities of $h$ has Lebesgue measure zero;
(A.3) for each fixed $\delta > 0$, $h$ is bounded on $[\delta, \pi]$;
(A.4) there exists $\alpha < 1$ such that $h(x) = O(|x|^{-\alpha})$ as $x \to 0$.

Given $h, k : I \to \mathbb{R}$, the *product of $h$ and $k$* is defined to be the function $(h \otimes k)(x,y) = h(x) \cdot k(y)$ for each $x, y \in I$.

LEMMA 3.2. *Assume that $f_s = f_{s1} \otimes f_{s2}$, $g_t = g_{t1} \otimes g_{t2}$, $f = \sum_{s=1}^m f_s$, $g = \sum_{t=1}^n g_t$ and each $f_{sj}, g_{tj}$ is admissible on $I, j = 1, 2$. Moreover, suppose that $f_{sj}(x) = O(|x|^{-\alpha})$, $g_{tj}(x) = O(|x|^{-\beta})$ as $x \to 0$ and $\alpha + \beta < \frac{1}{2}$. Then:*

(a) $\frac{1}{N^2} \mathrm{Tr}[T_N(f)T_N(g)]^2 \to (4\pi^2)^3 \int_{I^2} [fg]^2 \, dx\, dy,$
(b) $\frac{1}{N^k} \mathrm{Tr}[T_N(f)T_N(g)]^k \to 0$ *when* $k = 3, 4, \ldots.$

PROOF. First, observe that if $h_{st} = f_{s1} \otimes g_{t2}$, then $T_N(h_{st}) = T_N(f_{s1}) \otimes T_N(g_{t2}) = (a_{k\ell}T_N(g_{t2}))$. Indeed, a typical element in $(a_{k\ell}T_N(g_{t2}))$ is of the form $a_{k\ell}b_{mn} = \int_I e^{i(k-\ell)x} f_{s1}(x)\,dx \cdot \int_I e^{i(m-n)y} g_{t2}(y)\,dy = \int_{I^2} e^{i[(k-\ell)x+(m-n)y]} \times h_{st}(x,y)\,dx\,dy$. The latter quantity appears in row $N(k-1)+m$ and column $N(\ell-1)+n$ of both $T_N(f_{s1}) \otimes T_N(g_{t2})$ and $T_N(h_{st})$, which establishes the equality.

(a) Employing Theorem 1(a) of [12] and properties of the Kronecker product,

$$\frac{1}{N^2}\operatorname{Tr}[T_N(f)T_N(g)]^2$$

$$= \frac{1}{N^2}\sum_{s,t,u,v}\operatorname{Tr}[T_N(f_s)T_N(g_t)T_N(f_u)T_N(g_v)]$$

$$= \frac{1}{N^2}\sum_{s,t,u,v}\operatorname{Tr}[(T_N(f_{s1})\otimes T_N(f_{s2}))(T_N(g_{t1})\otimes T_N(g_{t2}))$$

$$\times (T_N(f_{u1})\otimes T_N(f_{u2}))\otimes(T_N(g_{v1})\otimes T_N(g_{v2}))]$$

$$= \sum_{s,t,u,v}\frac{1}{N}\operatorname{Tr}[T_N(f_{s1})T_N(g_{t1})T_N(f_{u1})T_N(g_{v1})]$$

$$\times \frac{1}{N}\operatorname{Tr}[T_N(f_{s2})T_N(g_{t2})T_N(f_{u2})T_N(g_{v2})]$$

$$\to \sum_{s,t,u,v}(2\pi)^3\int_I[f_{s1}g_{t1}f_{u1}g_{v1}]\,dx \cdot (2\pi)^3\int_I[f_{s2}g_{t2}f_{u2}g_{v2}]\,dy$$

$$= (4\pi^2)^3\int_{I^2}[fg]^2\,dx\,dy.$$

The above application of Theorem 1 of [12] is valid since each $f_{sj}(g_{tj})$ has the same order as $x \to 0$, respectively. It is not necessary that all $f_{sj}$'s be equal in the proof of Theorem 1 of [12].

(b) Suppose that $k$ is an integer exceeding 2. First, assume that $k(\alpha+\beta) < 1$. An extension of the argument used in part (a) shows that $\frac{1}{N^2}\operatorname{Tr}[T_N(f)T_N(g)]^k \to (4\pi^2)^{2k-1}\int_{I^2}[fg]^k\,dx\,dy$, and since $k > 2$, $\frac{1}{N^k}\operatorname{Tr}[T_N(f)T_N(g)]^k \to 0$. The case when $k(\alpha+\beta) \geq 1$ is verified in a similar manner by employing Theorem 1(b) of [12]. □

Given the process $\{X_{st} : s, t \in \mathbb{Z}\}$, recall that $\overline{X}_N = \frac{1}{N^2}\sum_{s,t=1}^N X_{st}$, and define $\tilde{X}'_N = (X_{11} - \overline{X}_N, X_{12} - \overline{X}_N, \ldots, X_{1N} - \overline{X}_N, X_{21} - \overline{X}_N, \ldots, X_{2N} - \overline{X}_N, \ldots, X_{N1} - \overline{X}_N, \ldots, X_{NN} - \overline{X}_N)$. Verification of the following result is based on the fact that the normal distribution is determined by its moments. The details of the proof are given in [6].

LEMMA 3.3. *Let $\{X_{st}: s,t \in \mathbb{Z}\}$ be a stationary Gaussian process having mean $\mu$ and spectral density function $f$, and let $g:(I^2, \mathfrak{B}_2) \to (\mathbb{R}, \mathfrak{B})$ be a bounded measurable function obeying $g(-x,-y) = g(x,y)$ for each $(x,y) \in I^2$. Moreover, assume that $f = \sum_{s=1}^{m} f_s$ $(g = \sum_{t=1}^{n} g_t)$, $f_s = f_{s1} \otimes f_{s2}$ $(g_t = g_{t1} \otimes g_{t2})$ and each $f_{sj}(g_{tj})$ is an admissible function with parameter $\alpha(\beta)$ in Definition* 3.1, *where $\alpha + \beta < \frac{1}{2}$. If $A_N = T_N(g)$, then $\frac{1}{N}[\tilde{X}_N' A_N \tilde{X}_N - E(\tilde{X}_N' A_N \tilde{X}_N)] \xrightarrow{D} \mathcal{N}(0, \delta^2)$, where $\delta^2 = 128\pi^6 \int_{I^2} [f(x,y) \cdot g(x,y)]^2 \, dx \, dy$.*

Recall that $\hat{\theta}_N$ is a value of $\theta = (\alpha, \beta, d) = (\alpha, \beta, d_1, d_2)$ making $\sigma_N^2(\theta)$ in (1.7) a minimum relative to the spectral density function $f(x,y,\theta)$ given in (1.6). For simplicity we use $\theta = (\theta_1, \theta_2, \theta_3, \theta_4) = (\alpha, \beta, d_1, d_2)$. Define a symmetric matrix

$$\Sigma = (\sigma_{ij})$$

(3.1)
$$\text{with } \sigma_{ij} = \int_{I^2} \frac{\partial f^{-1}(x,y,\theta_0)}{\partial \theta_i} \cdot \frac{\partial f^{-1}(x,y,\theta_0)}{\partial \theta_j} f^2(x,y,\theta_0) \, dx \, dy.$$

A calculation shows that $\sigma_{ij} = \int_{I^2} \frac{\partial^2 f^{-1}(x,y,\theta_0)}{\partial \theta_i \, \partial \theta_j} \cdot f(x,y,\theta_0) \, dx \, dy$ and, moreover, using the spectral density function listed in (1.6), it is straightforward to verify

(3.2)
$$\begin{aligned}
\sigma_{11} &= 8\pi^2/(1 - \alpha_0^2), & \sigma_{12} &= 0, \\
\sigma_{13} &= -8\pi^2(\text{Log}(1 - \alpha_0))/\alpha_0, & \sigma_{14} &= 0, \\
\sigma_{22} &= 8\pi^2/(1 - \beta_0^2), & \sigma_{23} &= 0, \\
\sigma_{24} &= -8\pi^2(\text{Log}(1 - \beta_0))/\beta_0, \\
\sigma_{33} &= \sigma_{44} = 4\pi^4/3 \quad \text{and} \quad \sigma_{34} = 0.
\end{aligned}$$

The interpretation given when $\alpha_0 = 0$ ($\beta_0 = 0$) is that $\sigma_{13}$ ($\sigma_{24}$) be replaced by its limiting value $8\pi^2$. Moreover, denote

$$C(x) = \frac{1}{1-x^2} \frac{\pi^2}{6} - \frac{\text{Log}^2(1-x)}{x^2}$$

and define

$$K = \frac{1}{8\pi^2 C(\alpha_0) C(\beta_0)},$$

$$\delta_{11} = \frac{\pi^2}{6} C(\beta_0), \qquad \delta_{12} = 0,$$

(3.3) $\quad \delta_{13} = \dfrac{\text{Log}(1-\alpha_0)}{\alpha_0} C(\beta_0), \qquad \delta_{14} = 0,$

FRACTIONAL SPATIAL PROCESSES 9

$$\delta_{22} = \frac{\pi^2}{6}C(\alpha_0), \qquad \delta_{23} = 0, \qquad \delta_{24} = \frac{\text{Log}(1-\beta_0)}{\beta_0}C(\alpha_0),$$

$$\delta_{33} = \frac{1}{1-\alpha_0^2}C(\beta_0), \qquad \delta_{34} = 0 \quad \text{and} \quad \delta_{44} = \frac{1}{1-\beta_0^2}C(\alpha_0).$$

Then one can verify that $\Sigma^{-1} = K(\delta_{ij})$.

PROOF OF THEOREM 1.2(a). Applying the mean value theorem,

(3.4) $$\frac{\partial}{\partial \theta}\sigma_N^2(\hat{\theta}_N) - \frac{\partial}{\partial \theta}\sigma_N^2(\theta_0) = \frac{\partial^2}{\partial \theta^2}\sigma_N^2(\overline{\theta}_N)(\hat{\theta}_N - \theta_0)$$

for some $\overline{\theta}_N$ belonging to the line segment between $\hat{\theta}_N$ and $\theta_0$. According to Theorem 1.1, $\hat{\theta}_N \to \theta_0$ almost surely and, thus, in probability; hence, $\frac{\partial}{\partial \theta}\sigma_N^2(\hat{\theta}_N) = 0$ since $\hat{\theta}_N$ becomes an interior point of the parameter space when $N$ is sufficiently large. It follows from Lemma 2.3, Theorem 1.1 and (3.1)–(3.2) that

(3.5) $$A_N := \frac{\partial^2}{\partial \theta^2}\sigma_N^2(\overline{\theta}_N) \to \Sigma \qquad \text{almost surely as } N \to \infty.$$

The Cramér–Wold device can be used to show that, as $N \to \infty$,

(3.6) $$N\left[\frac{\partial}{\partial \theta}\sigma_N^2(\theta_0) - E\left(\frac{\partial}{\partial \theta}\sigma_N^2(\theta_0)\right)\right] \xrightarrow{D} \mathcal{N}(\mathbf{0}, 128\pi^6\Sigma).$$

Indeed, let $c' = (c_1, c_2, c_3) \in \mathbb{R}^3$ and define $g(x,y,c) = c' \cdot \frac{\partial f^{-1}(x,y,\theta_0)}{\partial \theta}$; then

$$c' \cdot N\frac{\partial}{\partial \theta}\sigma_N^2(\theta_0) = N\int_{I^2} g(x,y,c)I_N(x,y)\,dx\,dy = \frac{1}{N}\tilde{X}_N' B_N \tilde{X}_N,$$

where $B_N = T_N(g)$. According to Lemma 3.3,

$$c'N\left[\frac{\partial}{\partial \theta}\sigma_N^2(\theta_0) - E\left(\frac{\partial}{\partial \theta}\sigma_N^2(\theta_0)\right)\right]$$
$$= \frac{1}{N}[\tilde{X}_N' B_N \tilde{X}_N - E(\tilde{X}_N' B_N \tilde{X}_N)] \xrightarrow{D} \mathcal{N}(0,\delta^2),$$

where

$$\delta^2 = 128\pi^6 \int_{I^2} [f(x,y,\theta_0) \cdot g(x,y,c)]^2\,dx\,dy$$
$$= 128\pi^6 \sum_{k,\ell=1}^3 c_k c_\ell \int_{I^2} \frac{\partial f^{-1}(x,y,\theta_0)}{\partial \theta_k} \cdot \frac{\partial f^{-1}(x,y,\theta_0)}{\partial \theta_\ell} f^2(x,y,\theta_0)\,dx\,dy$$
$$= 128\pi^6 c'\Sigma c.$$

Therefore, (3.6) holds. From this, employing (3.4) and (3.5),

$$N(\hat{\theta}_N - \theta_0) + NA_N^{-1}E\left(\frac{\partial}{\partial\theta}\sigma_N^2(\theta_0)\right) = -NA_N^{-1}\left[\frac{\partial}{\partial\theta}\sigma_N^2(\theta_0) - E\left(\frac{\partial}{\partial\theta}\sigma_N^2(\theta_0)\right)\right]$$

$$\xrightarrow{D} \mathcal{N}(\mathbf{0}, 128\pi^6\Sigma^{-1})$$

as $N \to \infty$. □

**4. Proof of Theorem 1.2(b).** According to (3.5), $A_N^{-1} \to \Sigma^{-1}$ almost surely and, thus, the proof of Theorem 1.2(b) involves determining the orders of $NE(\frac{\partial}{\partial\theta}\sigma_N^2(\theta_0))$. Define $\sigma_N^{2*}(\theta)$ as in (1.7) with $\overline{X}_N$ in $I_N$ replaced by $\mu = E(X_{11})$. The argument given by Fox and Taqqu ([12], Lemma 8.1) verifies that $NE[\frac{\partial}{\partial\theta}\sigma_N^2(\theta_0) - \frac{\partial}{\partial\theta}\sigma_N^{2*}(\theta_0)] \to 0$ as $N \to \infty$; it is shown that the orders of $NE(\frac{\partial}{\partial\theta}\sigma_N^{2*}(\theta_0))$ and $NE(\frac{\partial}{\partial\theta}\sigma_N^2(\theta_0))$ coincide.

Assume that $\{X_{st} : s, t \in \mathbb{Z}\}$ is a process described in Theorem 1.2 with spectral density function $f(x, y, \theta)$ defined in (1.6) and having mean $\mu$ and autocovariance function $\gamma(k, \ell)$. Denote

$$I_N^*(x, y) = \frac{1}{N^2}\left|\sum_{k,\ell=1}^N e^{i(kx+\ell y)}(X_{k\ell} - \mu)\right|^2$$

and

$$\sigma_N^{2*}(\theta) = \int_{I^2} \frac{I_N^*(x, y)}{f(x, y, \theta)}\, dx\, dy.$$

It follows from the stationarity of the $X$-process that

$$E\left(\frac{\partial}{\partial\theta}\sigma_N^{2*}(\theta_0)\right) = \frac{1}{N^2}\int_{I^2}\sum_{k,\ell,m,n=1}^N e^{i[(k-m)x+(\ell-n)y]}\frac{\partial}{\partial\theta}f^{-1}(x, y, \theta_0)$$

$$\times E(X_{k\ell} - \mu)(X_{mn} - \mu)\, dx\, dy$$

$$= \frac{1}{N^2}\sum_{k,\ell,m,n=1}^N \gamma(k-m, \ell-n)$$

(4.1)
$$\times \int_{I^2} e^{i[(k-m)x+(\ell-n)y]}\frac{\partial}{\partial\theta}f^{-1}(x, y, \theta_0)\, dx\, dy$$

$$= \frac{1}{N^2}\sum_{|k|<N, |\ell|<N}(N - |k|)(N - |\ell|)\gamma(k, \ell)$$

$$\times \int_{I^2} e^{i(kx+\ell y)}\frac{\partial}{\partial\theta}f^{-1}(x, y, \theta_0)\, dx\, dy.$$

Recall that $\theta'_0 = (\alpha_0, \beta_0, d_{10}, d_{20})$ denotes the true parameters. The following notation is used:

$$g(x, \alpha, d_1) = |1 - e^{-ix}|^{-2d_1}|1 - \alpha e^{-ix}|^{-2},$$
$$h(y, \beta, d_2) = |1 - e^{-iy}|^{-2d_2}|1 - \beta e^{-iy}|^{-2},$$

(4.2)
$$a_{k1} = \int_I e^{ikx} g(x, \alpha_0, d_{10})\, dx, \qquad a_{\ell 2} = \int_I e^{i\ell y} h(y, \beta_0, d_{20})\, dy,$$
$$b_{\ell 2} = \int_I e^{i\ell y} h^{-1}(y, \beta_0, d_{20})\, dy, \qquad c_{k1} = \int_I e^{ikx} \frac{\partial g^{-1}}{\partial d_1}(x, \alpha_0, d_{10})\, dx,$$

$$\psi(y) = \int_I g(y - x, \alpha_0, d_{10}) \frac{\partial g^{-1}}{\partial d_1}(x, \alpha_0, d_{10})\, dx.$$

It remains to determine the limit of the sequence $\{N(\text{Log }N)^{-3} E(\frac{\partial}{\partial \theta} \sigma_N^{2*}(\theta_0))\}$. For the sake of brevity, verification is presented here only for the $d_1$-component of $\theta$. Employing the notation in (4.2), $f^{-1}(x, y, \theta) = \frac{4\pi^2}{\sigma^2} g^{-1}(x, \alpha, d_1) \cdot h^{-1}(y, \beta, d_2)$ and, thus, (4.1) becomes

(4.3) $$E\left(\frac{\partial}{\partial d_1} \sigma_N^{2*}(\theta_0)\right) = \sum_{|k|<N}\left(1 - \frac{|k|}{N}\right) a_{k1} c_{k1} \cdot \sum_{|\ell|<N}\left(1 - \frac{|\ell|}{N}\right) a_{\ell 2} b_{\ell 2}.$$

Observe that the second summation converges to $8\pi^3$ as $N \to \infty$ since the limit is precisely $4\pi^2$ times the Cesàro sum of the convolution of $h$ and $h^{-1}$ evaluated at zero.

LEMMA 4.1. *Suppose that the assumptions listed in Theorem 2.1 are satisfied. Employing the notation of* (4.2) *gives the following:*

(i) $\psi(y) = -8d_{10} \sin \frac{y}{2} \text{Log}^2 \sin y + O(y|\text{Log }y|)$ as $y \to 0+$;

(ii) $E(\frac{\partial}{\partial \alpha} \sigma_N^2(\theta_0)) = O(\frac{\text{Log}^2 N}{N}), E(\frac{\partial}{\partial \beta} \sigma_N^2(\theta_0)) = O(\frac{\text{Log}^2 N}{N})$ and $\frac{N}{\text{Log}^3 N} \times E(\frac{\partial}{\partial d_i} \sigma_N^2(\theta_0)) \to -\frac{64}{3}\pi^2 d_{i0}, i = 1, 2,$ as $N \to \infty$.

PROOF. Verification of (i) involves applications of Lemmas A.1 and A.2 in the Appendix, together with several technical arguments of approximation and expansion. The complete details of the proof are available on request from the authors.

(ii). Proof of the third part for $i = 1$ is supplied here. According to (4.3), it remains to determine the limit of sequence $\frac{N}{\text{Log}^3 N} V_N := \frac{N}{\text{Log}^3 N} \Sigma_{|k|<N}(1 - \frac{|k|}{N}) a_{k1} c_{k1}$. Then $V_N$ is the $N$th Cesàro sum for $\psi$ evaluated at zero and thus

has the integral representation $V_N = \int_I \psi(y) K_N(y)\, dy$, where

$$(4.4) \qquad K_N(y) = \begin{cases} \dfrac{1}{2\pi N} \dfrac{\sin^2 Ny/2}{\sin^2 y/2}, & y \neq 0, \\ \dfrac{N}{2\pi}, & y = 0, \end{cases}$$

denotes the Fejér kernel (e.g., [7], page 71, or [26], page 88). Since $\psi$ is an even function, $\frac{1}{2} V_N = \int_{(0,1/N)} \psi(y) K_N(y)\, dy + \int_{[1/N,c)} \psi(y) K_N(y)\, dy + \int_{[c,\pi]} \psi(y) K_N(y)\, dy := J_{N1} + J_{N2} + J_{N3}$, where $c \in (0,1)$. Choosing $c$ sufficiently small, according to Lemma A.2, $\psi(y) = O(y \operatorname{Log}^2 y)$ as $y \to 0+$ and, thus, $|J_{N1}| \leq K_1 N \int_{(0,1/N)} y \operatorname{Log}^2 y\, dy = O(\frac{\operatorname{Log}^2 N}{N})$. Let $\|\psi\|_1$ denote the $L^1$-norm of $\psi$. Then $|J_{N3}| \leq \frac{K_3}{N} \|\psi\|_1 = O(\frac{1}{N})$.

It remains to estimate $J_{N2}$. Using the expansion for $\psi$ in (i) when $c$ is sufficiently small and $N$ sufficiently large,

$$J_{N2} = \int_{[1/N,c)} \left( -8 d_{10} \sin \frac{y}{2} \operatorname{Log}^2 \sin y + O(y|\operatorname{Log} y|) \right) K_N(y)\, dy.$$

The integral of the second term is $O(\frac{\operatorname{Log}^2 N}{N})$ and

$$\int_{[1/N,c)} \operatorname{Log}^2 \sin y \frac{\sin^2 Ny/2}{\sin y/2}\, dy$$

$$= \int_{[1/N,c)} \left( \frac{\operatorname{Log}^2 \sin y}{\sin y/2} - 2\frac{\operatorname{Log}^2 y}{y} \right) \sin^2 \frac{Ny}{2}\, dy$$

$$\quad + 2 \int_{[1/N,c)} \operatorname{Log}^2 y \frac{\sin^2 Ny/2}{y}\, dy$$

$$= \int_{[1/N,c)} O(y \operatorname{Log}^2 y) \sin^2 \frac{Ny}{2}\, dy + 2 \int_{[1/N,c)} \operatorname{Log}^2 y \frac{\sin^2 Ny/2}{y}\, dy$$

$$=: I_1 + I_2.$$

The first integral $I_1$ is $O(1)$ and $I_2 = \int_{[1/N,c)} \frac{\operatorname{Log}^2 y}{y}\, dy - \int_{[1/N,c)} \frac{\operatorname{Log}^2 y}{y} \cos Ny\, dy$. Integrating by parts, it is shown that the second integral is $O(\operatorname{Log}^2 N)$ and, thus, $I_2 = \frac{1}{3} \operatorname{Log}^3 N + O(\operatorname{Log}^2 N)$. Hence,

$$\int_{[1/N,c)} \operatorname{Log}^2 \sin y \frac{\sin^2 Ny/2}{\sin y/2}\, dy = I_1 + I_2 = O(1) + \frac{1}{3} \operatorname{Log}^3 N + O(\operatorname{Log}^2 N)$$

$$= \frac{1}{3} \operatorname{Log}^3 N + O(\operatorname{Log}^2 N)$$

and, thus,

$$J_{N2} = -\frac{8 d_{10}}{2\pi N} \int_{[1/N,c)} \operatorname{Log}^2 \sin y \frac{\sin^2 Ny/2}{\sin y/2}\, dy + O\left( \frac{\operatorname{Log}^2 N}{N} \right)$$

$$= -\frac{4d_{10}\operatorname{Log}^3 N}{3\pi N} + O\Big(\frac{\operatorname{Log}^2 N}{N}\Big).$$

Combining the above results,

$$V_N = 2(J_{N1} + J_{N2} + J_{N3})$$
$$= O\Big(\frac{\operatorname{Log}^2 N}{N}\Big) + \Big(-\frac{8d_{10}\operatorname{Log}^3 N}{3\pi N} + O\Big(\frac{\operatorname{Log}^2 N}{N}\Big)\Big) + O\Big(\frac{1}{N}\Big)$$
$$= -\frac{8d_{10}\operatorname{Log}^3 N}{3\pi N} + O\Big(\frac{\operatorname{Log}^2 N}{N}\Big)$$

and, thus, $\frac{N}{\operatorname{Log}^3 N}V_N \to -\frac{8d_{10}}{3\pi}$ as $N \to \infty$. It follows from (4.3) that

$$\frac{N}{\operatorname{Log}^3 N}E\Big(\frac{\partial}{\partial d_1}\sigma_N^{2\,*}(\theta_0)\Big) \to -\frac{8d_{10}}{3\pi}8\pi^3 = -\frac{64\pi^2 d_{10}}{3}.$$

As mentioned earlier, $NE(\frac{\partial}{\partial \theta}\sigma_N^2(\theta_0) - \frac{\partial}{\partial \theta}\sigma_N^{2\,*}(\theta_0)) \to 0$ and, hence,

$$\frac{N}{\operatorname{Log}^3 N}E\Big(\frac{\partial}{\partial d_1}\sigma_N^2(\theta_0)\Big) \to -\frac{64\pi^2 d_{10}}{3}. \qquad \square$$

PROOF OF THEOREM 1.2(b). Recall from (3.3) and (3.5) that $A_N^{-1} \to \Sigma^{-1} = K \cdot (\delta_{ij})$ almost surely. According to Theorem 1.2(a) and Lemma 4.1,

$$\frac{N}{\operatorname{Log}^3 N}(\hat{\theta}_N - \theta_0) + \frac{N}{\operatorname{Log}^3 N}A_N^{-1}E\Big(\frac{\partial}{\partial \theta}\sigma_N^2(\theta_0)\Big) \xrightarrow{P} 0$$

and $N(\operatorname{Log}^3 N)^{-1}E(\frac{\partial}{\partial \theta}\sigma_N^2(\theta_0)) \xrightarrow{P} \xi$, where $\xi' = -\frac{64\pi^2}{3}(0, 0, d_{10}, d_{20})$. It follows that

$$\frac{N}{\operatorname{Log}^3 N}(\hat{\theta}_N - \theta_0) \xrightarrow{P} -\Sigma^{-1} \cdot \xi.$$

In particular, $N(\operatorname{Log}^3 N)^{-1}(\hat{\alpha}_N - \alpha_0) \xrightarrow{P} -K \cdot \sum_{j=1}^4 \delta_{1j}\xi_j = -K\delta_{13}\xi_3$ since $\delta_{14} = 0$. Likewise, $N(\operatorname{Log}^3 N)^{-1}(\hat{\beta}_N - \beta_0) \xrightarrow{P} -K\delta_{24}\xi_4$ since $\delta_{23} = 0$. Further, $N(\operatorname{Log}^3 N)^{-1}(\hat{d}_1 - d_{10}) \xrightarrow{P} -K\delta_{33}\xi_3$, $N(\operatorname{Log}^3 N)^{-1}(\hat{d}_2 - d_{20}) \xrightarrow{P} -K\delta_{44}\xi_4$ and, hence, Theorem 1.2(b) follows. $\square$

## APPENDIX

The reader is referred to Boissy [6] for detailed proofs of the following lemmas.

LEMMA A.1. *Assume that $0 < d < 1$ is fixed. Then for $y > 0$ sufficiently small,*

(i) $0 \le d(\cos\frac{y}{2})^{1-d}\cot(\frac{x-y}{2})\sin\frac{y}{2} \le (\frac{\sin x/2}{\sin(x-y)/2})^d - (\cos\frac{y}{2})^d \le d(\cos\frac{y}{2})^{d-1} \times \cot(\frac{x-y}{2})\sin\frac{y}{2}$ when $x \in [-\pi, -\pi+y]$;

(ii) $0 \le -d\cot(\frac{x-y}{2})\sin\frac{y}{2} \le (\cos\frac{y}{2})^d - (\frac{\sin x/2}{\sin(x-y)/2})^d \le -2d\cot(\frac{x-y}{2})\sin\frac{y}{2}$ when $x \in [-\pi+y, -y]$;

(iii) $0 \le d2^{d-1}\cot(\frac{x-y}{2})\sin\frac{y}{2} \le (\frac{\sin x/2}{\sin(x-y)/2})^d - (\cos\frac{y}{2})^d \le d(\cos\frac{y}{2})^{d-1} \times \cot(\frac{x-y}{2})\sin\frac{y}{2}$ when $x \in [2y, \pi]$.

LEMMA A.2. *Fix $0 < d < 1$. Then there exist positive constants $C_1, C_2$ and $\delta$ such that, for $0 < y < \delta$, $C_1 y(\text{Log } y)^2 \le |\int_I |2\sin\frac{x-y}{2}|^{-d}|2\sin\frac{x}{2}|^d \times \text{Log} |2\sin\frac{x}{2}| dx| \le C_2 y(\text{Log } y)^2$.*

**Acknowledgments.** The authors express their indebtedness to the referee and an Associate Editor for comments leading to significant improvements in the manuscript.

Y. BOISSY
X. LI
G. D. RICHARDSON
DEPARTMENT OF MATHEMATICS
UNIVERSITY OF CENTRAL FLORIDA
ORLANDO, FLORIDA 32816-1364
USA
E-MAIL: xli@pegasus.cc.ucf.edu
    garyr@pegasus.cc.ucf.edu

B. B. BHATTACHARYYA
DEPARTMENT OF STATISTICS
NORTH CAROLINA STATE UNIVERSITY
RALEIGH, NORTH CAROLINA 27695
USA
E-MAIL: bhattach@stat.ncsu.edu